\documentclass{amsart}
\usepackage[english]{babel}
\usepackage{amsmath,amsthm}
\input xy
\xyoption{all}
\newcommand{\norm}[1]{\left\Vert#1\right\Vert}
\newcommand{\R}{\mathbb{R}}
\newcommand{\lie}[1]{\mathfrak{#1}}     %Lie algebras
\newcommand{\Lie}{\mathcal{L}}          %Lie derivative
\newcommand{\Z}{\mathbb{Z}}
\newcommand{\C}{\mathbb{C}}
\newcommand{\hook}{\lrcorner\,}
\newcommand{\LieG}[1]{\mathrm{#1}}      %Lie groups
\newcommand{\SU}{\mathrm{SU}}
\newcommand{\SO}{\mathrm{SO}}
\newcommand{\Gtwo}{\mathrm{G}_2}
\newcommand{\GL}{\mathrm{GL}}
\newcommand{\SL}{\mathrm{SL}}
\newcommand{\dfn}[1]{\emph{#1}}
\newcommand{\leftdoublebracket}{[\![}
\newcommand{\rightdoublebracket}{]\!]}
\DeclareMathOperator{\Hom}{Hom}
\DeclareMathOperator{\End}{End}
\DeclareMathOperator{\ad}{ad}
\renewcommand{\Re}{\mathfrak{Re}\,}
\renewcommand{\Im}{\mathfrak{Im}\,}

\theoremstyle{plain}
\newtheorem{proposition}{Proposition}[section]
\newtheorem{theorem}[proposition]{Theorem}
\newtheorem{lemma}[proposition]{Lemma}

\theoremstyle{definition}
\newtheorem{definition}[proposition]{Definition}

\theoremstyle{remark}
\newtheorem*{remark}{Remark}
\begin{document}
\title{Special symplectic six-manifolds}
%\date{\today}
\author{Diego Conti}
\address{Scuola Normale Superiore\\Piazza dei Cavalieri 7
\\ 56126 Pisa\\ Italy}
\email{diego.conti@sns.it}
\author{Adriano Tomassini}
\address{Dipartimento di Matematica\\ Universit\`a di Parma\\ Parco Area delle Scienze 53/A\\
43100 Parma\\ Italy} \email{adriano.tomassini@unipr.it}
\begin{abstract} We classify  nilmanifolds with an invariant symplectic half-flat structure. We solve the half-flat
evolution equations in one  example,  writing down the resulting Ricci-flat metric. We study the geometry of the orbit
space of 6-manifolds with an $\SU(3)$-structure preserved by a $\mathrm{U}(1)$ action, giving  characterizations in the
symplectic half-flat and integrable case.
\end{abstract}
\maketitle
\section{Introduction}
\label{sec:introduction} \sloppy A half-flat manifold is a six-dimensional manifold endowed with an
$\SU(3)$\nobreakdash-struc\-ture whose intrinsic torsion is symmetric. An $\SU(3)$-structure defines a
non\nobreakdash-\hspace{0pt}degenerate two-form $\omega$, an almost-complex structure $J$, and a complex volume form
$\Psi$; the half-flat condition is equivalent to requiring $\omega\wedge\omega$ and the real part of $\Psi$ to be
closed \cite{ChiossiSalamon}. \\
Hypersurfaces in seven-dimensional manifolds with holonomy $\Gtwo$ have a natural
half-flat structure, given by the restriction of the holonomy group representation; the intrinsic torsion can then be
identified with the second fundamental form. The converse is not obvious. In \cite{Hitchin:StableForms}, Hitchin proved
that, starting with a half-flat manifold $(M,\omega,\Psi)$, if certain evolution equations admit a solution coinciding
with $(\omega,\Psi)$ at time zero, then $(M,\omega,\Psi)$ can be embedded isometrically as a hypersurface in a manifold
with holonomy contained in $\Gtwo$.\\
 Given a six-manifold $M$ with an $\SU(3)$-structure, one can consider the product $\Gtwo$-structure on $M\times S^1$.
Half-flat manifolds satisfying special conditions on this product $\Gtwo$-structure have been studied: in
\cite{ChiossiSwann} and \cite{ChiossiFino},  the six-dimensional nilmanifolds carrying invariant $\SU(3)$-structures of
these types have been classified. More generally, the  problem of classifying nilmanifolds admitting invariant
half-flat structures is  open.\\
In this paper we focus on the symplectic case; that is to say, we take into consideration half-flat structures for
which $\omega$ is closed. In this context, one can introduce special Lagrangian submanifolds, namely the
three-dimensional submanifolds $\iota\colon L\to M$ such that both  $\iota^*\omega$ and $\iota^*\Im\Psi$ vanish. Like
in the integrable case, special Lagrangian submanifolds are exactly the submanifolds calibrated by $\Re\Psi$
\cite{HarveyLawson:CalibratedGeometries}. Symplectic half-flat manifolds can also be viewed as symplectic manifolds
$(M,\omega)$ endowed with two extra objects, namely an $\omega$-calibrated almost-complex structure and a $(3,0)$-form
with closed real part which is parallel with respect to
the Chern connection \cite{dBT3}, \cite{dBT2}.\\
Our main result is the classification of nilmanifolds carrying an invariant symplectic half-flat structure; the special
case $b_1\geq 4$ was carried out in \cite{Bedulli}.

\fussy The contents of this paper are organized as follows. Section~\ref{sec:InvariantGCY} consists of the
classification. Since six-dimensional nilmanifolds can be realized as circle bundles over nilmanifolds of dimension
five, this classification problem can be reduced to a problem in five dimensions: indeed, the symplectic half-flat
structure induces an $\SU(2)$-structure on the base of the circle bundle, satisfying certain conditions involving the
curvature (see Lemma~\ref{lemma:GCY}). In Lemma~\ref{lemma:GCY2} we classify the five-dimensional nilmanifolds with
this type of induced structure. We then use this lemma to show that every six-dimensional nilmanifold admitting a
symplectic half-flat structure is modelled on one of a list of three Lie algebras (Theorem~\ref{thm:InvariantGCY}). The
corresponding nilmanifolds are the torus, a torus bundle over the four-dimensional torus (see \cite{dBT2} and
\cite{Giovannini}) and a torus bundle over a three-dimensional torus (see
\cite{Bedulli} and \cite{Giovannini}). In the last case, the fibres are actually special Lagrangian submanifolds.\\
In Section~\ref{sec:Evolution} we fix a symplectic half-flat structure on the irreducible nilmanifold appearing in
Theorem~\ref{thm:InvariantGCY}, and solve the evolution equations. We write down the resulting metric, computing the
curvature and proving that the holonomy coincides with $\Gtwo$. We also observe that if one starts with the reducible
example, where the half-flat manifold is  the product of a five-manifold with a circle, the
resulting seven-dimensional manifold is also reducible.\\
In Section~\ref{sec:Functor} we generalize this situation and the construction of Lemma~\ref{lemma:GCY} to the
non-invariant case. More precisely, we consider circle bundles with a $\LieG{U}(1)$\nobreakdash-invariant
$\SU(3)$-structure, and we define an induced $\SU(2)$-structure on the base; the manifolds are assumed to be compact.
The metric underlying the $\SU(2)$-structure is not the same as the quotient metric,
but it is obtained from it by rescaling along certain directions, as in \cite{ApostolovSalamon}.\\
The compactness assumption is dropped in Section~\ref{sec:Intrinsictorsion}, where we compute the intrinsic torsion of
the $\SU(2)$-structure in terms of the intrinsic torsion of the $\SU(3)$\nobreakdash-structure and  curvature of the
bundle. As far as we know, this is the first detailed application of intrinsic torsion for $\SU(2)$-structures on
five-manifolds. We then characterize the $\LieG{U}(1)$\nobreakdash-\hspace{0pt}invariant symplectic half-flat manifolds
in terms of the quotient structure (Proposition~\ref{prop:GCYFromHypo}). In the case that the $\LieG{U}(1)$-invariant
$\SU(3)$\nobreakdash-structure is integrable, we obtain a stronger result: indeed, the intrinsic torsion of the
quotient structure and the curvature
 are determined by the length of the fibres (Theorem~\ref{thm:QuotientOfCY}).
\section{Invariant structures on nilmanifolds}
\label{sec:InvariantGCY} In this section we introduce half-flat structures on $6$-manifolds and classify
invariant symplectic half-flat structures on 6-dimensional nilmanifolds.\newline
Let $M$ be a 6-dimensional manifold. An $\SU(3)$-structure on $M$ is a pair
$(\omega,\Psi)$, where $\omega$ is a non-degenerate two-form and $\Psi$ is a decomposable complex three-form,
such that the following compatibility conditions hold: \begin{equation} \label{eqn:CompatibilitySU3}\left\{
\begin{aligned}
\Psi\wedge\omega&=0\\
\Psi\wedge\overline\Psi&=\frac43 i\,\omega^3
\end{aligned}\right.
\end{equation}
Indeed, a decomposable complex three-form $\Psi=\theta^1\wedge\theta^2\wedge\theta^3$ determines an almost-complex
structure $J$ for which $\theta^1$, $\theta^2$ and $\theta^3$ span the space of forms of type $(1,0)$. The second
compatibility condition implies in particular that $J$ is $\omega$-tamed; the first  condition asserts
that $\omega$ is of type $(1,1)$, so $J$ is actually calibrated by $\omega$.\\
At each point $x$ of $M$, the three-form $\Psi_x$ has stabilizer conjugate to $\SL(3,\C)$ for the natural
action of $\GL^+(T_xM)$; on the other hand, the stabilizer of $\omega_x$ is conjugate to the symplectic group
$\LieG{Sp}(3,\R)$ for the natural action of $\GL(T_xM)$. The compatibility conditions ensure that the
intersection of the two stabilizers is conjugate to $\SU(3)$, so that an $\SU(3)$-structure is defined.\\
We shall denote by $\psi^+$ and $\psi^-$ the real and imaginary part of $\Psi$, respectively. It was shown by Hitchin
\cite{Hitchin:TheGeometryOfThreeForms} that having fixed the orientation, the real three-form $\psi^+$ is sufficient to
determine the almost-complex structure $J$, and therefore  $\psi^-$. Thus, an $\SU(3)$-structure $(\omega,\Psi)$ is
really determined by the pair $(\omega,\psi^+)$.

We now introduce a special class of $\SU(3)$-structures on 6-manifolds,  related to 7-dimensional Riemannian manifolds
with holonomy contained in $\Gtwo$ \cite{ChiossiSalamon}:
\begin{definition}
An $\SU(3)$-structure $(\omega,\psi^+)$ on a 6-manifold is \dfn{half-flat} if $\omega\wedge\omega$ and
$\psi^+$ are closed.
\end{definition}
The 2-form $\omega$ appearing in the characterization of $\SU(3)$-structures is required to be non-degenerate; if it
is also closed, it defines a symplectic structure. In this case, we say that the $\SU(3)$-structure is symplectic.

Consider a nilmanifold $M$, i.e. a compact manifold of the form $\Gamma\backslash G$, where $G$ is a 6-dimensional
nilpotent group, and $\Gamma$ a discrete subgroup of $G$. Recall that in six dimensions, every nilpotent
Lie algebra $\lie{g}$ gives rise  to such a nilmanifold.\\
We say that a structure on the nilmanifold $M$ is \dfn{invariant} if it pulls back to a left-invariant structure on
$G$. Invariant structures can be viewed as structures on the Lie algebra $\lie{g}$; using the above characterization of
$\SU(3)$-structures in terms of differential forms, we shall mainly work with the
dual $\lie{g}^*$.\\
We start by reducing the problem to a problem in five dimensions; the idea is to realize $M$ as a circle bundle over a
5-dimensional manifold, in such a way that the $\SU(3)$-structure on $M$ is invariant under the circle action. The
geometry of this construction will be studied in Section~\ref{sec:Functor}; here, we shall use the following algebraic
result:
\begin{lemma}\label{lemma:GCY}
Let $(\omega,\psi^+)$ be a symplectic half-flat structure on  a nilpotent Lie algebra $\lie{g}$; we have an orthogonal
decomposition
\[\lie{g}^*=\langle \eta\rangle\oplus V^5 \;,\]
where $\eta$ is a unit form, and
\[d (\lie{g}^*)\subseteq \Lambda^2V^5\;.\] Define forms $\alpha$, $\omega_1$, $\omega_2$ and $\omega_3$ on
$\ker\eta$ by
\begin{equation}\label{eqn:HypoProduct}
\left\{\begin{aligned}
    \omega&=\omega_3+\eta\wedge\alpha\\
    \Psi&=(\omega_1+i\omega_2)\wedge (\eta+i\alpha)
\end{aligned}\right.
\end{equation}
Setting $\phi=d\eta$, the following hold:
\begin{equation}\label{eqn:GCYFromHypoUnitXi}\left\{
\begin{aligned}
d\alpha&=0\;, & d\omega_1&=0\;,\\
d\omega_3&=-\phi\wedge\alpha\;,&    d(\omega_2\wedge\alpha)&=\omega_1\wedge\phi\;.
\end{aligned}\right.
\end{equation}
\end{lemma}
\begin{proof}
By Engel's theorem, some non-zero $\xi$ in $\lie{g}$ satisfies \[\ad(\xi)=0\;.\] Choosing for $\eta$ a suitable
multiple of $\xi^\flat$, and setting $V^5=\eta^\perp$, the first part of the Lemma is satisfied.\\
 By definition,
\[0=d\omega= d\omega_3+\phi\wedge\alpha-\eta\wedge d\alpha\;;\]
isolating the component in $\eta\wedge \Lambda^2 V^5$, we deduce that $\alpha$ is closed and $d\omega_3$
satisfies the required equation.\\
Similarly, the rest of \eqref{eqn:GCYFromHypoUnitXi} follow from:
\[0=d\psi^+=d\omega_1\wedge\eta+\omega_1\wedge\phi - d(\omega_2\wedge\alpha)\;.\qedhere\]
\end{proof}
\begin{remark}
The forms $(\alpha,\omega_i)$ introduced in Lemma~\ref{lemma:GCY} define an $\SU(2)$\nobreakdash-\hspace{0pt}structure
on $\lie{g}$. More generally, we recall that differential forms $(\alpha,\omega_1,\omega_2,\omega_3)$ on a 5-manifold
define an $\SU(2)$-structure if and only if at each point (and hence locally) there exists a coframe $e^1,\dots,e^5$
such that
\begin{equation}
\label{eqn:referenceformsSU2}\left\{
\begin{aligned}\alpha&=e^5&\omega_1&=e^{12}+e^{34}\\
\omega_2&=e^{13}+e^{42}&\omega_3&=e^{14}+e^{23}
\end{aligned}\right.
\end{equation}
Here and in the sequel, $e^{12}$ is short for $e^1\wedge e^2$, and so on. If one fixes an orientation, the condition
above is equivalent to the existence of a  triplet $(\omega_1,\psi_2,\psi_3)$ with
\begin{align*}\omega_1&=e^{12}+e^{34}\;, & \psi_2&=e^{135}+e^{425}\;,&\psi_3&=e^{145}+e^{235}\;.
\end{align*}
\end{remark}
\smallskip
By construction, $V^5$ is itself the dual of a nilpotent Lie algebra; we shall proceed by listing the 5-dimensional Lie
algebras that arise this way.  To describe Lie algebras, we shall use symbolic expressions such as
\[\lie{g}=(0,0,0,0,0,12)\;,\] meaning that $\lie{g}^*$ has a basis $\eta^1,\dots,\eta^6$ such that
$d\eta^6=\eta^1\wedge\eta^2$ and $\eta^i$ is closed for $i=1,\dotsc,5$.
\begin{lemma}\label{lemma:GCY2}
In the hypotheses of Lemma \ref{lemma:GCY}, $V^5$ is one of
\[(0,0,0,0,0)\;,\quad (0,0,0,0,12)\;, \quad (0,0,0,12,13)\;.\]
\end{lemma}
\begin{proof}
By the same argument we used in the proof of Lemma~\ref{lemma:GCY}, we can construct a filtration
\[V^0\subset\dotsb\subset
V^5\;,\quad \dim V^i=i\;,\;  dV^{i+1}\subset\Lambda^2 V^i\;.\]
Moreover, by Lemma~\ref{lemma:GCY}, we can assume that
$\alpha$ lies in $V^1\subset V^4$; therefore, using the fact that the 4-dimensional representation of $\SU(2)$ is
transitive, we can fix a basis $e^1,\dotsc,e^5$ of $V^5$ satisfying \eqref{eqn:referenceformsSU2}, with
$e^4$ in $(V^4)^\perp$.\\
We have to show that the first Betti number satisfies $b_1\geq 3$; by the classification of 5-dimensional nilpotent Lie
algebras, it will then suffice to show that for some choice of the $V^i$ as above, $de^4$ lies in $\Lambda^2V^3$ (in
particular, this implies that the step of $V^5$ is less than two).

 By Lemma~\ref{lemma:GCY},
\[0=de^{12}+de^{34} \equiv de^3\wedge e^4 \mod \Lambda^2V^4\;,\]
implying that $e^3$ is closed. Thus, we can assume
\[V^2=\langle e^3,e^5\rangle\;;\quad V^4=\langle e^1,e^2,e^3,e^5\rangle\;.\]
Define a real constant $h$, a 2-form $\gamma\in\Lambda^2\langle e^1,e^2,e^3\rangle$ and 1-forms $\phi_4$,~$\phi_5$ in
$\langle e^1,e^2,e^3\rangle$, such that
 \[\phi=\phi_5\wedge e^5+  \phi_4\wedge e^4 + he^{45} +
\gamma\;.\]
 Since $\phi$ is closed,
 \begin{equation}\label{eqn:GCY:dphi}d\phi_5\wedge e^5  - \phi_4\wedge de^4 + hde^4\wedge e^5 + d\gamma=- d\phi_4\wedge e^4
 \;;\end{equation}
the left-hand side lies in $\Lambda^3 V^4$, i.e. it has no component containing $e^4$, so both sides are zero and
$d\phi_4=0$.\\
By Lemma~\ref{lemma:GCY} $d\omega_3=-\phi\wedge\alpha$, giving
\begin{align}
\label{eqn:GCY2}de^{14}+de^{2}\wedge e^3=-\phi_4\wedge e^{45} - \gamma\wedge e^5\;;
\end{align}
comparing the components in $e^4\wedge\Lambda^2V^4$, we obtain
\[de^1= \phi_4\wedge e^5\;.\]
Again by Lemma~\ref{lemma:GCY}, $d\psi_2=\omega_1\wedge\phi$; on the other hand, $d\psi_2=d\omega_2\wedge\alpha$, so we
can drop the component of $\phi$ not containing $e^5$ and write
\begin{align}
 \label{eqn:GCY3}(de^{13}+de^{42})\wedge
e^5&=(e^{12}+e^{34})\wedge(\phi_5\wedge e^5 + he^{45})\;.
\end{align}
The components  containing $e^4$ give
\[de^2\wedge e^5=e^3\wedge\phi_5\wedge e^5-  he^{125}\;,\]
and wedging with $e^3$,
\[de^2\wedge e^{35}=-h e^{1235}\;.\]
Since $de^2$ is in $\Lambda^2V^3$, the left hand side is in $\Lambda^4V^3$, so it is zero. We conclude that $h=0$  and
 \begin{equation}
 \label{eqn:GCY3b}
de^2\wedge e^5=e^3\wedge\phi_5\wedge e^5\;;
\end{equation}
now, either $V^3=\langle e^1,e^3,e^5\rangle$, or some linear combination $\lambda e^1+ e^2$ lies in $V^3$, and
consequently $0=(\lambda de^1+de^2)\wedge e^5=de^2\wedge e^5\;.$ Either way, \[\phi_5\in\langle e^1,e^3\rangle\;.\] By
Lemma~\ref{lemma:GCY} $\omega_1$ is closed, giving
 \begin{equation} \label{eqn:GCY1b}
 \phi_4\wedge e^{52}-e^1\wedge de^2=e^3\wedge de^4\;;
\end{equation}
to proceed further, we must distinguish three cases.\\

\emph{i}) Suppose that $\phi_4$ is not a multiple of $e^3$; then \[V^3=\langle e^3,e^5,\phi_4\rangle\;,\] and $d$ is
zero on $V^3$. Moreover $e^1$ is not closed, so $V^4=V^3\oplus \langle e^1\rangle$. Since $e^2$ is in $V^4$, we have
$de^2=k\,e^5\wedge\phi_4$ for some (possibly zero) constant $k$. So \eqref{eqn:GCY1b} becomes
\[\phi_4\wedge e^{52}-k\,e^{15}\wedge\phi_4=e^3\wedge de^4\;,\]
implying that \[de^4\wedge e^3\wedge \phi_4=0=de^4\wedge e^{35}\;.\] Since the space of closed two-forms in
$\Lambda^2V^4$ is
\[\Lambda^2V^3\oplus \langle e^{15}, e^1\wedge\gamma_4\rangle\;,\]
we can conclude that $de^4$ lies in $\Lambda^2V^3$; we already know that $b_1\geq 3$, so there is nothing left to prove
in this case.\\

In general, the component of \eqref{eqn:GCY2} in $\Lambda^3V^4$ gives
\begin{equation}
\label{eqn:GCY2b} -e^1\wedge de^4+de^2\wedge e^3=-\gamma\wedge e^5\;;
\end{equation}
in particular, $de^4\wedge e^{15}=0$. Moreover, we can rewrite \eqref{eqn:GCY:dphi} as
 \begin{equation}
 \label{eqn:GCY4b}-\phi_4\wedge de^4+d\gamma=0\;.\end{equation}

\emph{ii}) \sloppy Suppose now that $\phi_4=0$; then $e^1$ is closed and we can assume that \mbox{$V^3=\langle
e^1,e^3,e^5\rangle$}. By \eqref{eqn:GCY1b},
\[-e^1\wedge de^2=e^3\wedge de^4\;,\] \fussy
so clearly $de^4\wedge e^{13}=0$; moreover $de^4\wedge e^{35}=0$, since by \eqref{eqn:GCY3b}  $de^2\wedge e^{15}$ is
zero.
It follows that $de^4$ lies in $\Lambda^2V^3$, completing the proof in this case.\\

\emph{iii}) The remaining case  is the one where $\phi_4=a e^3$ for some non-zero $a$. By \eqref{eqn:GCY4b}
and~\eqref{eqn:GCY1b}, this condition implies
\begin{equation}\label{eqn:dgammatilde}
d\gamma=ae^3\wedge de^4=a^2e^{352}-ae^1\wedge de^2\;.
\end{equation}
Equations \eqref{eqn:GCY3b} and \eqref{eqn:GCY1b} show that
 \[de^2\wedge e^{15}=0=de^2\wedge e^{13}=de^2\wedge e^{35}\;,\]
so $de^2$ lies in $\Lambda^2\langle e^1,e^3,e^5\rangle$. Hence the space of closed forms in $\Lambda^2V^4$ is contained
in
\[\Lambda^2\langle e^1,e^3,e^5\rangle\oplus e^2\wedge\langle e^3,e^5\rangle\;;\]
wedging the closed two-form $\gamma-ae^{12}$ with $e^{35}$ we deduce \[\gamma\wedge e^{35}=ae^{1235}\;.\] Comparing
with \eqref{eqn:GCY2b}, we find $e^{13}\wedge de^4=-ae^{1235}$, which together with \eqref{eqn:dgammatilde} gives the
contradiction
\[-a^2e^{1235}=a^2e^{1352}\;.\qedhere\]
\end{proof}

All three possibilities listed in Lemma~\ref{lemma:GCY2} can occur:
\begin{itemize}
\item On $V^5=(0,0,0,0,0)$, set
 \[\omega_1=\eta^{12}+\eta^{34}\;,\quad \psi_2=\eta^{135}+\eta^{425}\;,\quad
 \psi_3=\eta^{145}+\eta^{235}\;,\quad \phi=0\;.\]
\item
On $V^5=(0,0,0,0,12)$, set
 \[ \omega_1=\eta^{34}+\eta^{15}\;,\quad \psi_2=\eta^{312}+\eta^{542}\;,\quad
 \psi_3=\eta^{352}+\eta^{412}\;, \quad\phi=-\eta^{13}\;.\]
\item
On $V^5=(0,0,0,12,13)$, set \[ \begin{aligned}
 \omega_1&=\eta^{24}+\eta^{35}\;,& \psi_2&=\eta^{123}+\eta^{154}\;,&
 \psi_3&=\eta^{125}+\eta^{143}\;,&\phi&=-2\eta^{23}\;; \text{ or}\\
\omega_1&=\eta^{24}-\eta^{35}\;,& \psi_2&=-\eta^{123}+\eta^{154}\;,&
 \psi_3&=\eta^{125}-\eta^{143}\;,&\phi&=0\;.
\end{aligned}\]
\end{itemize}
It is easy to verify that Equations \ref{eqn:GCYFromHypoUnitXi} are satisfied in these cases. The construction can then
be inverted: define
\[\lie{g}^*=\langle\eta\rangle\oplus V^5\;,\] declaring that $d\eta=\phi$; clearly, $\lie{g}^*$ is the dual of a nilpotent Lie algebra $\lie{g}$. A
straightforward calculation shows that the $\SU(3)$-structure on $\lie{g}$ defined by \eqref{eqn:HypoProduct} is
half-flat and symplectic.

So, there are  three non-isomorphic nilpotent Lie algebras that admit a symplectic half-flat structure.  It only
remains to show that this list is complete.
\begin{theorem}
\label{thm:InvariantGCY} The 6-dimensional nilpotent Lie algebras whose corresponding nilmanifold carries an invariant
symplectic half-flat structure are
\[
 (0,0,0,0,0,0)\;,\quad
 (0,0,0,0,12,13)\;,\quad(0,0,0,12,13,23)\;.
\]
\end{theorem}

\begin{proof}
We retain the  notation from  the proof of Lemma~\ref{lemma:GCY2}.  We first show that  $\phi_4$ is zero; in other
words, case \emph{i}) of  Lemma~\ref{lemma:GCY2} cannot occur, like case  \emph{iii}), which we have already ruled out.
Indeed, suppose that $\phi_4$ is independent of $e^3$ and $e^5$. Hence $V^5=(0,0,0,12,13)$; indeed, $de^1$ and $de^4$
must be independent, since otherwise a combination of $e^1$ and $e^4$ would lie in $\ker d$, which is orthogonal
to $e^4$ by construction.\\
Observe that $\langle e^1,e^3,\phi^4\rangle$ has dimension three, because $\phi_4$ is closed but $e^1$ is not, and we
are assuming that $\phi_4$ is not a multiple of $e^3$. Therefore,
\[\langle e^1,e^3,\phi^4\rangle=\langle e^1,e^2,e^3\rangle\;,\]
and we can write $\gamma=a\,e^{13}+b\,e^1\wedge\phi_4+c\,e^3\wedge\phi_4$; Equation \ref{eqn:GCY2b} then yields
\[-e^1\wedge de^4+ke^5\wedge \phi_4\wedge e^3=-a\,e^{135}+b\,e^{15}\wedge\phi_4+c\,e^{35}\wedge\phi_4\;,\]
so in particular
\[de^4=a\,e^{35}-b\,e^5\wedge \phi_4\;;\]
substituting in \eqref{eqn:GCY4b},  it follows that
\[-a\,\phi_4\wedge e^{35}+a\, \phi_4\wedge e^{53}=0\;,\]
i.e. $a=0$; but then $de^4=b\,e^5\wedge\phi_4$ is a multiple of $de^1$, which is absurd.

We have proved that $\phi_4$ is necessarily zero; now assume that $e^2$ is closed. Then \eqref{eqn:GCY1b} and
\eqref{eqn:GCY2b} give
\[ e^3\wedge de^4=0\;,\quad e^1\wedge de^4=\gamma\wedge e^5\;.\]
It follows that $de^4=\lambda\, e^{35}+ \mu\, e^{13}$ and $\gamma = \lambda\, e^{13}$ for some constants $\lambda$ and
$\mu$. The components of  \eqref{eqn:GCY3} not containing $e^4$ give
 \[\mu\, e^{1325} = -e^{125}\wedge\phi_5\;;\]
Since as a consequence of  \eqref{eqn:GCY3b} $\phi_5$ is a multiple of $e^3$, it follows that $\phi_5=-\mu\, e^3$.
Summing up,
\[\phi=-\mu\, e^{35} +\lambda\, e^{13}\;,\]
so that $\phi$ and $de^4$ are either linearly independent or both zero. The resulting 6-dimensional Lie algebras are
\[(0,0,0,0,12,13)\;,\quad (0,0,0,0,0,0)  \;.\]

If $e^2$ is not closed, $V^5$ is $(0,0,0,12,13)$. As both $de^4$ and $de^2$ are in $\Lambda^2V^3$,
\eqref{eqn:GCY2b}~implies that $\gamma$ is a multiple of $e^{13}$. Therefore $\phi$ lies in $\Lambda^2V^3$ as well,
forcing $\lie{g}$ to be either $(0,0,0,12,13,23)$ or $(0,0,0,0,12,13)$.
\end{proof}
\section{Associated Ricci-flat metrics}
\label{sec:Evolution} In this section we show that the symplectic half-flat manifolds of Theorem \ref{thm:InvariantGCY}
can be realized as hypersurfaces in Ricci-flat manifolds; in one example, we compute explicitly the metric and its
curvature, proving that the holonomy is $\Gtwo$.\\
Recall that a $\Gtwo$-structure on a 7-manifold  $N$ is defined by a three-form $\varphi$ which at each point $x$ lies
in the $\GL(T_xN)$ orbit of
 \begin{equation}
 \label{eqn:G2Form}
 e^{147}+e^{257}+e^{367}+e^{123}-e^{156}-e^{426}-e^{453}\;,
  \end{equation}
where $e^1,\dotsc,e^7$ is any coframe at $x$. Since $\Gtwo$ is contained in $\SO(7)$, $\varphi$ determines a metric and
an orientation on $N$.\\
The $\Gtwo$-structure defined by $\varphi$ is integrable if and only if $\varphi$ is closed and co-closed
\cite{FernandezGray}. Then the corresponding Riemannian metric  has  holonomy contained in $\Gtwo$ and is therefore
Ricci-flat. \\
Now let $N$ be a manifold  with an integrable $\Gtwo$-structure and let  $\iota\colon M\to N$ be a hypersurface. Then
there exists
 a unique $\SU(3)$-structure $(\omega,\psi^+)$ on $M$ such that
\[\psi^+=\iota^*\varphi\;,\quad \omega^2=2\iota^*\ast\!\varphi\;;\]
from  the integrability of the $\Gtwo$-structure it clearly follows  that this induced structure is half-flat.
Conversely, it is known that if $(\omega(t),\psi^+(t))$ is a one-parameter family of half-flat structures on $M$, for
$t$ ranging in $(a,b)$, then
\begin{equation}
\varphi=\omega\wedge dt+\psi^+
\end{equation}
defines a $\Gtwo$-structure on $M\times(a,b)$, which is integrable if and only if $(\omega(t),\psi^+(t))$  satisfies
the \dfn{half-flat evolution equations}:
 \begin{equation}\label{eqn:HalfFlatEvolution}
 d\omega=\frac\partial{\partial t}\psi^+\;,\quad d\psi^-=-\frac{1}{2}\frac\partial{\partial t}\omega^2\;.
 \end{equation}
In this construction, $(\omega(t),\psi^+(t))$ must satisfy the compatibility conditions \eqref{eqn:CompatibilitySU3}
for all $t$. However, it turns out that $\psi^-$ is still defined for small deformations of an $\SU(3)$-structure, and
if  \eqref{eqn:HalfFlatEvolution} are satisfied, then \eqref{eqn:CompatibilitySU3} are preserved in time. Indeed, we
have the following  \cite{Hitchin:StableForms}:
 \begin{theorem}[Hitchin]
 \label{thm:HitchinEvolution}
Let $M$ be a compact 6-manifold, and let  $(\omega(t),\psi^+(t))$ be a one-parameter family of sections of
$\Lambda^2(M)\oplus\Lambda^3(M)$ satisfying the evolution equations \eqref{eqn:HalfFlatEvolution}. If
$(\omega(0),\psi^+(0))$ is a half-flat structure, and $\omega(t)^3$ is nowhere zero for $t$ in $(a,b)$, then
$(\omega(t),\psi^+(t))$ defines a half-flat structure for all $t\in(a,b)$. In particular, $M\times(a,b)$ has a
Riemannian metric with holonomy contained in $\Gtwo$.
\end{theorem}
We now solve the evolution equations \eqref{eqn:HalfFlatEvolution} for the nilmanifold with Lie algebra
$(0,0,0,12,13,23)$. Our solution is different from the one given in \cite{ChiossiFino}, since we choose symplectic
initial data. Consider the one-parameter family of $\SU(3)$-structures given by
\begin{equation} \label{eqn:EvolSol1}\left\{
\begin{aligned}
  \omega&=\frac{1}{u}\eta^{16}-\frac{1}{u}\eta^{25}-\frac{3u^2-1}u \eta^{34}\\
  \psi^+&=\frac{(3u^2-1)^2}{4u^6} \eta^{123}-2\eta^{154}+2\eta^{624}+\eta^{653}
\end{aligned}\right.
\end{equation}
Clearly, $(\omega,\psi^+)$ is half-flat for all values of $u$, and symplectic for $u=\pm 1$. Setting
\[t=-12+\frac1{2u^3}-\frac1{10u^5}\;,\]
Equations \ref{eqn:EvolSol1} give a solution of \eqref{eqn:HalfFlatEvolution}. An orthonormal basis of 1-forms is given
by
\begin{align*}
E^1&=\sqrt\frac{3u^2-1}{2u^4} \eta^1\;,& E^2&=\sqrt\frac{3u^2-1}{2u^4}\eta^2\;,& E^3&=\frac{3u^2-1}{2u^2}\eta^3\;, \\
E^4&=\sqrt\frac{2u^2}{3u^2-1} \eta^6\;,& E^5&=-\sqrt\frac{2u^2}{3u^2-1} \eta^5\;, &E^6&=-2u\,\eta^4\;,
\end{align*}
where indices and signs have been adjusted for compatibility with \eqref{eqn:G2Form}. Let $N$ be the 7-dimensional
Riemannian manifold obtained by the above construction; set $E^7=dt$, and consider the inclusion
$\Lambda^2(N)\subset\End(TN)$, where the two-form $E^{ij}$ is identified with the skew-symmetric endomorphism mapping
$E^i$ to $E^j$.  As a section of $S^2(\Lambda^2(N))$, the curvature is given by
\begin{multline*}
 -\frac{4u^{10}}{(3u^2-1)^4}\left( 3\left(E^{17}+E^{35}\right)^2 + 3\left(E^{34}-E^{27}\right)^2 +\left(E^{14}-E^{25}\right)^2 - \left(E^{12}+E^{45}\right)
 \right)+\\
 -\frac{12u^{10}(2u^2-1)}{(3u^2-1)^3}\left( \left(E^{16}+E^{27}\right)^2 + \left(E^{17}-E^{26}\right)^2
 -2\left(E^{12}-E^{67}\right)^2 \right)+\\
 +\frac{12u^{10}(u^2-1)}{(3u^2-1)^4}\left(\left(E^{24}+E^{37}\right)^2+\left(E^{15}-E^{37}\right)^2\right)+\\
 -\frac{12u^{10}(u^2-2)}{(3u^2-1)^4}\left(\left(E^{13}+E^{57}\right)^2 + \left(E^{23}-E^{47}\right)^2\right)+\\
 -\frac{4u^{10}}{(3u^2-1)^3}\left(\left(E^{23}+E^{56}\right)^2 + \left(E^{13}+E^{46}\right)^2 -
 \left(E^{14}-E^{36}\right)^2\right)\;.
\end{multline*}
This shows that the metric is not reducible and the holonomy equals $\Gtwo$.\\

Similar computations for a symplectic half-flat structure on the Lie algebra \mbox{$\lie{g}=(0,0,0,0,12,13)$} were
carried out in \cite{ContiSalamon}. In this case though, as one can see by computing the curvature, the resulting
7-manifold is the product of a 6-manifold with holonomy $\SU(3)$ and a circle. In fact, the nilmanifold is a trivial
circle bundle with connection form $\eta =\eta^4$. The symplectic half-flat structure induces an
$\SU(2)$\nobreakdash-structure $(\alpha,\omega_i)$ on the five-dimensional base by \eqref{eqn:HypoProduct}; this
structure is \dfn{hypo} in the sense of \cite{ContiSalamon} (see also Section~\ref{sec:Functor}), and can therefore be
evolved to give a 6\nobreakdash-manifold with holonomy $\SU(3)$. The reducible 7-manifold is nothing but the product of
this 6-manifold with a circle.

\begin{remark}
Whilst by Hitchin's theorem the half-flat conditions are automatically preserved in time, the symplectic condition is
not, as shown in the above example. This is a general fact: the evolution flow is transverse to the space of
symplectic half-flat structures, except where it vanishes (namely, at points defining integrable structures).
\end{remark}
\section{$\SU(3)$-structures on  circle bundles}
\label{sec:Functor} In this section we pursue an idea introduced in Section~\ref{sec:InvariantGCY}, namely that of
reducing  a 6-dimensional  manifold to a 5-dimensional manifold  by means of a quotient, and establishing a relation
between the two geometries in terms of $G$-structures. Here we work in a more general context, without requiring
invariance under a transitive action; however, for the construction to make sense we still need invariance along one
direction, i.e. a Killing field on the 6-manifold. More precisely, we shall establish a one-to-one correspondence
between a class of 6-manifolds with an $\SU(3)$-structure and a regular vector field preserving the structure, and a
class of 5-manifolds with an $\SU(2)$\nobreakdash-\hspace{0pt}structure plus some additional data; since this
correspondence only holds ``up to isomorphism'', it will be natural to state it in terms of categories. For the moment,
we  impose no integrability conditions on the $\SU(3)$-structure.\\
 We define a category $\mathcal{K}$ whose objects are 4-tuples $(M,\omega,\psi^+,X)$, where $M$ is a
compact 6-dimensional manifold, $(\omega,\psi^+)$ is an $\SU(3)$-structure on $M$, and  $X$ is a regular vector field
on $M$ which preserves the $\SU(3)$\nobreakdash-structure, i.e. \[\Lie_X\omega=0=\Lie_X\psi^+\;.\] For brevity, we
shall often write $M$ for $(M,\omega,\psi^+,X)$; so, when $M$ is referred to as an object of $\mathcal{K}$, it will be
understood that $\omega$, $\psi^+$ and $X$ are also fixed on $M$. Sometimes we will need to consider two distinct
objects $M$,
$\tilde M$, and it will be understood that $\tilde M$ stands for  $(\tilde M,\tilde\omega,\tilde\psi^+,\tilde X)$.\\
A morphism $f\in\Hom(M,\tilde M)$ is a smooth map $f\colon M\to\tilde M$ such that
 \begin{align*}
 f^*\tilde\omega&=\omega\;,& f^*\tilde\psi^+&=\psi^+\;, & X \text{ is $f$-related to } \tilde X\;.
 \end{align*}
In particular, morphisms are orientation-preserving isometries, and therefore covering maps.

We are going to relate $\mathcal{K}$ to a category $\mathcal{C}$, whose objects are 6-tuples
\[(N,\alpha,\omega_1,\omega_2,\omega_3,\phi,t)\;,\] where $N$ is a compact 5-manifold, $(\alpha,\omega_i)$ is an
$\SU(2)$-structure on $N$, $t$ is a function on $N$, and $\phi$ is a closed two-form on $N$ such that \[
\left[\frac{1}{2\pi}\phi\right]\in H^2(N,\Z)\;.\] Again, we shall write $N$
for an object $(N,\alpha,\omega_1,\omega_2,\omega_3,\phi,t)$ of $\mathcal{C}$.\\
A morphism $f\in\Hom(N,\tilde N)$ is a smooth map $f\colon N\to\tilde N$ such that
\begin{align*}
f^*\tilde\alpha&=\alpha\;,& f^*\tilde\omega_i&=\omega_i\;,\; i=1,2,3& f^*\tilde\phi&=\phi\;,&\tilde t\circ f&=t\;.
\end{align*}

We shall construct a functor $F\colon\mathcal{K}\to\mathcal{C}$ which realizes each object $M$ of $\mathcal{K}$ as a
circle bundle over $F(M)$; the 2-form $\phi$ represents the curvature, and the function $t$ the length of the fibres.\\
Let us first recall that a functor $F\colon\mathcal{C}\to\mathcal{D}$ is \dfn{faithful} (resp. \dfn{full})  if
\[F\colon \Hom(A,B) \to\Hom(F(A),F(B))\] is one-to-one (resp. onto) for all objects $A$, $B$ of $\mathcal{C}$; it is \dfn{representative}
if every object in $\mathcal{D}$ is isomorphic to $F(A)$ for some object $A$ of $\mathcal{C}$. A full, representative,
faithful functor is called an \dfn{equivalence}.\\
 The functor we  consider does not quite establish an
equivalence, but we shall show that it {induces} an equivalence between categories derived from $\mathcal{K}$ and
$\mathcal{C}$. The starting point is the following observation: if $M$ is an object of $\mathcal{K}$, the maximal
integral curves of $X$ are closed subsets of the compact manifold $M$, and therefore diffeomorphic to circles. More
precisely, we have the following \cite[p. 27]{Blair}:
\begin{lemma}\label{lemma:OrbitLength}
For every object $M$ of $\mathcal{K}$, the maximal integral curves of $X$ viewed as maps $\phi_x\colon \R\to M$ are
periodic, with period not depending on the point $x$. In particular $M$ is the total space of a circle bundle, and $X$
is a fundamental vector field.
\end{lemma}%
We can now prove the following:
\begin{proposition}\label{prop:Functor}
There is a  representative functor $F\colon\mathcal{K}\to\mathcal{C}$.
\end{proposition}
\begin{proof}
We define the covariant functor $F\colon\mathcal{K}\to\mathcal{C}$ as follows: let $M$ be an object of $\mathcal{K}$.
The space of integral lines of $X$ is a compact 5-manifold $N$, and by Lemma~\ref{lemma:OrbitLength} $M$ is the total
space of a circle bundle over $N$. Since $X$ is a Killing vector field, the norm of $X$ is constant on the fibres, and
so
defines a smooth function $t$ on $N$.\\
By Lemma~\ref{lemma:OrbitLength}, the maximal integral curves of $X$ have constant period; we can rescale $X$ so that
this period is $2\pi$. Let $\eta$ be the connection 1-form determined by $X$, i.e. the dual form to $X$ rescaled so
that $\eta(X)=1$; set $\phi=d\eta$. Then the cohomology class
\[c_1=\left[\frac{1}{2\pi} \phi\right]\]
is the Chern class of the $\LieG U(1)$-bundle $M\to N$; as such, it is integral. Define forms $(\alpha,\omega_i)$ on
$M$ by
\begin{equation}
\label{eqn:quotient2}\left\{
\begin{aligned}
\alpha&=X\hook\omega\;,  & \omega_3&=tX\hook(\omega\wedge\eta)\;,\\
\omega_1&=X\hook\psi^+\;, & \omega_2&=X\hook\psi^-\;.
\end{aligned}\right.
\end{equation}
%If one defines $\psi_2$ and $\psi_3$ in the usual way, it is clear that equations (\ref{eqn:quotient}) hold.
By construction, all the objects appearing on the right-hand sides of \eqref{eqn:quotient2} are invariant under the
action of $\LieG U(1)$; therefore, each form $\alpha$, $\omega_1$, $\omega_2$, $\omega_3$ is  the pullback of
a  form on $N$, which we denote by the same symbol.\\
Now choose a local orthonormal basis of 1-forms
 \[\frac 1{\sqrt t}e^1,\dots,\frac 1{\sqrt t}e^4,\frac1te^5,e^6\] on $M$, with $\eta=t^{-1}\,e^6$, such that
\begin{align*}
\omega&=\frac 1t (e^{14}+ e^{23}+e^{65})\;, & \Psi=\frac 1t(e^1+ie^4)\wedge(e^2+ie^3)\wedge(e^6+\frac 1ti e^5)\;.
\end{align*}
Then (\ref{eqn:referenceformsSU2}) is satisfied, so $(\alpha,\omega_i)$ defines an $\SU(2)$-structure on $N$. In
particular, $e^1,\dots,e^5$ is an orthonormal basis of 1-forms on $N$: so, we are not using the  quotient metric on the
5-manifold, but a deformation of it.

If $M$ and $\tilde M$ are objects of $\mathcal{K}$ and $f\in\Hom(M,\tilde M)$, then $f$ maps integral curves of $X$ to
integral curves of $\tilde X$. Therefore, $f$ induces a smooth map $F(f)\colon N\to \tilde N$, where $N=F(M)$, $\tilde
N=F(\tilde M)$; we must show that $F(f)$ is a morphism. From the fact that $X$ is $f$-related to  $\tilde X$ and $f$ is
a local isometry, it follows that $\tilde t\circ F(f)=t$. Now consider the diagram of maps
\[ \xymatrix{M\ar[d]^\pi\ar[r]^f & \tilde M\ar[d]^{\tilde\pi}\\ N\ar[r]^{F(f)}& \tilde N}\]
 By the commutativity of the diagram and \eqref{eqn:quotient2},
\[\pi^* \left( F(f)^*\tilde\alpha\right)=f^*(\tilde\pi^*\tilde\alpha)=\alpha\;;\]
 therefore, $F(f)^*\tilde\alpha=\alpha$ on $N$. By the same argument $F(f)^*\tilde\omega_i=\omega_i$ and
\mbox{ $F(f)^*\tilde\phi=\phi$}.\\

Next we show that $F$ is representative; more precisely, that every object $N$ of $\mathcal{C}$ can be written as
$F(M)$ for some object $M$ of $\mathcal{K}$. Indeed,  let $M$ be a circle bundle over $N$ with Chern class
$[\phi/(2\pi)]$. Let $A$ be the standard generator of the Lie algebra $\lie{u}(1)$, and let $X=A^*$ be the associated
fundamental vector field. Choose a connection form $\eta$ such that $d\eta=\phi$ and define
\begin{equation}\label{eqn:HypoProduct2}
\left\{\begin{aligned}
    \omega&=t^{-1}\omega_3+\eta\wedge\alpha\\
    \Psi&=(\omega_1+i\omega_2)\wedge (\eta+it^{-2}\alpha)
\end{aligned}\right.
\end{equation}
It is clear that this defines an $\SU(3)$-structure preserved by $X$, which has norm $t\eta(X)=t$, and Equations
\ref{eqn:quotient2} give back the original $\SU(2)$-structure on $N$.
\end{proof}
\begin{remark}
One might wonder at the advantage of referring explicitly to categories, rather than just defining $F$ and studying its
properties. With the latter approach, a problem arises when one tries to construct an inverse to $F$: indeed, ``the''
circle bundle with a given Chern class is not a well-defined manifold, but something only defined up to isomorphism.
\end{remark}
\begin{remark}
The reduction we have chosen behaves well with respect to evolution theory.  Indeed, let $N$ be a 5-manifold with an
$\SU(2)$-structure; recall that $N$ is called \dfn{hypo} if the forms $\omega_1$, $\omega_2\wedge\alpha$, and
$\omega_3\wedge\alpha$ are closed. It is easy to verify that the $\SU(3)$-structure induced on $N\times S^1$ by the
above construction (corresponding to taking $\phi=0$ and $t=1$) is half-flat if and only if the structure on $N$ is
hypo. Moreover,  hypo geometry also has evolution equations similar to \eqref{eqn:HalfFlatEvolution}, and it turns out
that a one-parameter solution of the hypo evolution equations lifts to a solution of the half-flat evolution equations.
An example of this situation is the reducible half-flat nilmanifold  mentioned in Section~\ref{sec:Evolution}.
\end{remark}

The functor $F$ fails to be an equivalence in two respects: it is not full, and it is  not faithful. We start by
addressing the first issue. Let $\mathcal{C}'$ be the subcategory of $\mathcal{C}$ consisting of objects $N$ such that
$N$ is simply-connected as a manifold; let $\mathcal{K}'$ be the subcategory of $\mathcal{K}$ of objects $M$ such that
$F(M)$ is simply-connected.
\begin{lemma}\label{lemma:full}
The functor $F\colon\mathcal{K}'\to\mathcal{C}'$ is full.
\end{lemma}
\begin{proof}
 Consider two objects $M$, $\tilde M$
in $\mathcal{K}'$; let $F(M)=N$, $F(\tilde M)=\tilde N$, and let $h\in \Hom(N,\tilde N)$. Fix a point $u$ in $M$; every
path $\sigma$ in $N$ based at $\pi(u)$ has a horizontal lift $\gamma$ with $\gamma(0)=u$. Fix also a point $\tilde u$
in $\tilde M$, lying over $\pi(u)$; then $\sigma$ lifts to $\tilde\gamma$ with $\tilde\gamma(0)=\tilde u$. Now define
$f(\gamma(t))=\tilde\gamma(t)$; we have to show that the definition does not
depend on the path $\sigma$.\\
Indeed, let $I=[0,1]$, and let $\Sigma\colon I\times I\to N$ be a smooth homotopy  with fixed endpoints between two
paths $\Sigma_0$ and $\Sigma_1$ starting at $\pi(u)$, where \mbox{$\Sigma_t=\Sigma(t,\cdot)$}. The pullback bundle
$\Sigma^*M$ is trivial; we can therefore choose a section  \[\Gamma\colon I\times I\to M\;,\] where $\Gamma_0$ and
$\Gamma_1$ are horizontal
lifts of $\Sigma_0$, $\Sigma_1$ respectively, and $\Gamma(\cdot,0)=u$.\\
The ``vertical distance'' from $\Gamma_0(1)$ to $\Gamma_1(1)$ is measured by
\[\int_{\Gamma(\cdot,1)} \eta = \int_{I\times I} \Gamma^*d\eta=\int_{I\times I} \Sigma^*\phi\;,\]
where we have used Stokes' theorem.
 Quite similarly if $\tilde\Sigma=h\circ\Sigma$ and $\tilde\Gamma$ is constructed
as above, imposing this time $\tilde\Gamma(\cdot,0)=\tilde u$, we obtain
\[\int_{\tilde\Gamma(\cdot,1)} \tilde\eta  =\int_{I\times I} \tilde\Gamma^*d\tilde\eta=\int_{I\times I} \tilde\Sigma^*\tilde\phi=\int_{I\times I}\Sigma^*\phi\;,\]
because $h^*\tilde\phi=\phi$.\\
 Hence, we can lift $h$ to an equivariant map $f\colon M\to\tilde M$ satisfying $f^*\tilde\eta=\eta$. Thus,
\[F\colon\Hom(M,\tilde M)\to \Hom(N,\tilde N)\] is onto as required.
\end{proof}
\begin{remark}
The definition of the lift $f$ in the proof of Lemma~\ref{lemma:full} is also a characterization, because morphisms
preserve the connection form, and therefore map horizontal paths to horizontal paths. So, in the non-simply-connected
case one cannot expect to be able to produce a lift, i.e. $F$ is not full as a functor from $\mathcal{K}$
to~$\mathcal{C}$.
\end{remark}

We now come to faithfulness. Fix an object $M$ in $\mathcal{K}'$ and consider the periodic integral lines $\phi_x$; the
map
\[ M\ni x \xrightarrow{f_\theta} \phi_x(\theta)\in M\]
is an isomorphism for all constants $\theta$. Clearly, $F(f_\theta)$ is the identity; so,  \[F\colon\Hom(M, M)\to
\Hom(F(M),F(M))\] is not one-to-one and $F$ is not faithful.\\
However,  this is the only amount to which $F$ fails to be an equivalence. Indeed, for each $M$, $\tilde M$ in
$\mathcal{K}'$, define an equivalence relation in $\Hom(M,\tilde M)$ by
\[f_1 \sim f_2 \iff f_1=f_2 \circ f_{\theta}\;,\;\theta\in\R\;.\]
We can then consider the \dfn{quotient category} $\mathcal{K}''$, whose objects are the objects of $\mathcal{K}'$, and
whose morphisms are defined by \[\Hom_{\mathcal{K}''}(M,\tilde M)=\Hom_{\mathcal{K}'}(M,\tilde M)/\sim\;.\]
\begin{proposition}
$\mathcal{K}''$ is a category, and the induced functor $F''\colon\mathcal{K}''\to \mathcal{C}'$ is an equivalence.
\end{proposition}
\begin{proof}
By construction, elements of $\Hom(M,\tilde M)$ are $\LieG U(1)$-equivariant. Therefore, \[(f_1 \circ
f_{\theta_1})\circ (f_2\circ f_{\theta_2})= f_1\circ f_2\circ f_{\theta_1+\theta_2}\;.\] It follows that if $g_1 \sim
f_1$ and $g_2\sim
f_2$ then $g_1\circ g_2\sim f_1\circ f_2$. This is sufficient to conclude that $\mathcal{K}''$ is a category.\\
Observe that the induced functor $F''$ is well defined, because clearly $f_1\sim f_2$ implies \mbox{$F''(f_1)=F''(f_2)$}.\\
Now, recall from the proof of Proposition~\ref{prop:Functor} that the connection form is defined only by the metric and
the Killing field; since a morphism \mbox{$f\in\Hom(M,\tilde M)$} is an isometry and $X$ is $f$-related to $\tilde X$,
we have $f^*\tilde\eta=\eta$. Therefore, $f$ is uniquely determined by its value at a point, or in other words,
\[F''\colon\Hom_{\mathcal{K}''}(M, M)\to \Hom(F''(M),F''(M))\] is  one-to-one, as required.
\end{proof}
\section{Intrinsic torsion of the quotient structure}
\label{sec:Intrinsictorsion} In this section we drop the assumptions of compactness and global regularity, and study
the local behaviour of the  construction of Section~\ref{sec:Functor} in terms of intrinsic torsion. In particular, we
characterize the intrinsic torsion of the $\SU(2)$-structures obtained by taking the quotient of a symplectic half-flat
structure, generalizing Lemma~\ref{lemma:GCY}. Then, in the assumption that the starting  $\SU(3)$-structure is
integrable, we write down a differential equation that the function $t$ must satisfy, and prove that the intrinsic
torsion of the quotient $\SU(2)$-structure depends only on $t$. More precisely, we give necessary and sufficient
conditions on $(N,\alpha,\omega_i,t)$ for it to arise, locally, as the quotient of a 6-manifold with an integrable
$\SU(3)$-structure. Observe that in the integrable case the 6-manifold cannot be compact, unless it is reducible.

We shall work in a neighbourhood of a point where the Killing field is non-zero; thus, we assume that $M$ is a
6-manifold with an $\SU(3)$-structure preserved by some regular Killing field $X$. Recall from the proof of Proposition
\ref{prop:Functor} that the quotient $N$ is a 5\nobreakdash-manifold on which an $\SU(2)$-structure $(\alpha,\omega_i)$
is induced (see \eqref{eqn:quotient2}), as well as a function $t$, the norm of $X$, and a two-form $\phi$, which in the
case of a circle
bundle is the curvature form.\\
Recall from \cite{ChiossiSalamon} that the intrinsic torsion of an $\SU(3)$-structure takes values in a 42-dimensional
space, and its components can be represented as follows:
\begin{equation}
\begin{array}{|c|c|}
\hline W_1^+&W_1^-\\
\hline\raisebox{0cm}[3.5mm]{}
W_2^+&W_2^-\\
\hline
\multicolumn{2}{|c|}{\raisebox{0cm}[3.5mm]{}W_3}\\
\hline
\multicolumn{2}{|c|}{\raisebox{0cm}[3.5mm]{}W_4}\\
\hline \multicolumn{2}{|c|}{\raisebox{0cm}[3.5mm]{}W_5}\\
\hline
\end{array}\quad\in\quad
\begin{array}{|c|c|}
\hline
\R&\R\\
\hline\raisebox{0cm}[3.5mm]{}
[\Lambda^{1,1}_0]&[\Lambda^{1,1}_0]\\
\hline
\multicolumn{2}{|c|}{\raisebox{0cm}[3.5mm]{}\leftdoublebracket\Lambda^{2,1}_0\rightdoublebracket}\\
\hline
\multicolumn{2}{|c|}{\raisebox{0cm}[3.5mm]{}\leftdoublebracket\Lambda^{1,0}\rightdoublebracket}\\
\hline \multicolumn{2}{|c|}{\raisebox{0cm}[3.5mm]{}\leftdoublebracket\Lambda^{1,0}\rightdoublebracket}\\
\hline
\end{array}
\end{equation}
meaning that the component $W_1^+$ takes values in $\R$, and so on. Explicitly, we can write
\begin{equation}\label{eqn:SU3Torsion}
\left\{\begin{aligned}
d\psi^+&=\psi^+\wedge W_5 \,+\, W_2^+\wedge\omega \,+\, W_1^+ \omega^2\\
d\psi^-&=\psi^-\wedge W_5 \,+\, W_2^-\wedge\omega \,+\, W_1^-  \omega^2\\
d\omega&=-\frac32 W_1^-\psi^+ \,+\, \frac32 W_1^+\psi^-\, +\, W_3\, +\, W_4\wedge\omega
\end{aligned}
\right.
\end{equation}
 We can do the same for $\SU(2)$-structures on 5-manifolds \cite{ContiSalamon}; the intrinsic
torsion now takes vales in a 35-dimensional space, and we can arrange its components in the following table:
\[
\begin{array}{|c|c|c|}
\hline  \multicolumn{3}{|c|}{\lambda}\\
\hline f_1 & f_2 & f_3\\
\hline\raisebox{0cm}[3.5mm]{} g_1^2 & g_1^3 & g_2^3\\
\hline \multicolumn{3}{|c|}{\beta}\\
\hline \gamma_1 & \gamma_2 & \gamma_3\\
\hline \multicolumn{3}{|c|}{\raisebox{0cm}[3mm]{}\omega^-}\\
\hline \sigma_1^- & \sigma_2^- & \sigma_3^- \\
\hline
\end{array}
\quad\in\quad
\begin{array}{|c|c|c|}
\hline  \multicolumn{3}{|c|}{\R}\\
\hline\R & \R & \R\\
\hline \R & \R & \R\\
\hline \multicolumn{3}{|c|}{\raisebox{0cm}[3.5mm]{}\Lambda^1}\\
\hline \raisebox{0cm}[3.5mm]{} \Lambda^1 & \Lambda^1 & \Lambda^1\\
\hline \multicolumn{3}{|c|}{\raisebox{0cm}[3.5mm]{}\Lambda^2_-}\\
\hline \raisebox{0cm}[3.5mm]{}\Lambda^2_- &\Lambda^2_- &\Lambda^2_-  \\
\hline
\end{array}
\]
In the above table, $\Lambda^1$ is the 4-dimensional representation of $\SU(2)$ such that the tangent space at a point
is
\[T=\Lambda^1\oplus\R\;,\] whereas $\Lambda^2_-$ is the 3-dimensional representation of $\SU(2)$ consisting of anti-self-dual
two-forms on  $\Lambda^1$. We shall write, say, $(\omega)_{\Lambda^2_-}$ for the ${\Lambda^2_-}$ component of a
two-form $\omega$.\\
Setting  $g_i^j=-g^i_j$, the components of the intrinsic torsion are  given by
 \begin{equation}
 \label{eqn:SU2Torsion}\left\{
\begin{aligned}
d\alpha&=\alpha\wedge\beta+\sum_{j=1}^3 f^j \omega_j+\omega^-\\
d\omega_i&=\gamma_i\wedge\omega_i+\lambda\,\alpha\wedge\omega_i+\sum_{j\neq i}
g_i^j\alpha\wedge\omega_j+\alpha\wedge\sigma_i^-
\end{aligned}\right.
\end{equation}
We can now prove the following:
\begin{proposition}\label{prop:QuotientTorsion}
Define the intrinsic torsion of $M$ as above, and write
 \begin{equation}
 \label{eqn:W_i}
 W_i = \eta\wedge \Xi_i + \Delta_i\;, \quad\text{where}\quad\Xi_i = X\hook W_i\;;
 \end{equation}
then the intrinsic torsion of the quotient is given by
\[
\begin{array}{|c|c|c|}
\hline  \multicolumn{3}{|c|}{\raisebox{0cm}[3.5mm]{} -\langle \Delta_5,\alpha\rangle}\\
\hline\raisebox{0cm}[3.5mm]{} \frac32W_1^- & -\frac32W_1^+ - \frac12\langle \Xi_3,\omega_2\rangle & -t^{-1}\Xi_4 -\frac12\langle \Xi_3,\omega_3\rangle\\
\hline\raisebox{0cm}[3.5mm]{} -t^{-2}\Xi_5 & -2t^{-1}W_1^+ -t^{-1}\langle \Xi_2^+,\alpha\rangle & -2t^{-1}W_1^- - t^{-1}\langle\Xi_2^-,\alpha\rangle\\
\hline \multicolumn{3}{|c|}{\raisebox{0cm}[3.5mm]{}-\bigl(\Delta_4\bigr)_{\Lambda^1} - \alpha\hook \Xi_3}\\
\hline \raisebox{0cm}[3.5mm][1.5mm]{}-\bigl(\Delta_5\bigr)_{\Lambda^1} -t^{-1}\Xi_2^+\hook \omega_2 & -\bigl(\Delta_5\bigr)_{\Lambda^1}+t^{-1}\Xi_2^-\hook \omega_1 & \bigl(\Delta_4+d\log t+\frac12t\omega_3\hook \Delta_3\bigr)_{\Lambda^1}\\
\hline \multicolumn{3}{|c|}{\raisebox{0cm}[3.5mm][1.5mm]{} -\bigl(\Xi_3\bigr)_{\Lambda^2_-}}\\
\hline \raisebox{0cm}[3.5mm]{} -\Delta_2^+ &-\Delta_2^- & t\bigl(\alpha\hook \Delta_3-\phi\bigr)_{\Lambda^2_-} \\
\hline
\end{array}\]
\end{proposition}
\begin{proof}
Taking the interior product of \eqref{eqn:SU3Torsion} with $X$, then substituting \eqref{eqn:quotient2} in the
left-hand side and \eqref{eqn:HypoProduct2} in the right-hand side, one easily computes:
\begin{align*}
d\alpha&=\frac32W_1^-\omega_1 - \frac32W_1^+\omega_2-\Xi_3-t^{-1}\Xi_4\, \omega_3+\Delta_4\wedge\alpha\\
d\omega_1&=-\omega_1\wedge \Delta_5-t^{-2}\Xi_5\,\omega_2\wedge\alpha -
t^{-1}\Xi_2^+\wedge\omega_3-\Delta_2^+\wedge\alpha-2t^{-1}W_1^+\omega_3\wedge\alpha \\
d\omega_2&=-\omega_2\wedge \Delta_5+t^{-2}\Xi_5\,\omega_1\wedge\alpha-
t^{-1}\Xi_2^-\wedge\omega_3-\Delta_2^-\wedge\alpha -2t^{-1}W_1^-\omega_3\wedge\alpha
\end{align*}
On the other hand,  $d\omega_3=d\log t\wedge\omega_3-tX\hook d(\omega\wedge\eta)$ by \eqref{eqn:quotient2}. Using
\eqref{eqn:SU3Torsion} and then \eqref{eqn:HypoProduct2}, we obtain
 \begin{multline*}
 d(\omega\wedge\eta)= \frac32t^{-2}W_1^-\omega_2\wedge\alpha\wedge\eta+\frac32t^{-2}
W_1^+\omega_1\wedge\alpha\wedge\eta+\Delta_3\wedge\eta\,+\\
+t^{-1}\Delta_4\wedge\omega_3\wedge\eta+t^{-1}\omega_3\wedge\phi+\eta\wedge\alpha\wedge\phi\;;
 \end{multline*}
therefore
\[ d\omega_3=d\log t\wedge\omega_3 +\frac32t^{-1}W_1^-\omega_2\wedge\alpha
+\frac32t^{-1}W_1^+\omega_1\wedge\alpha +t\Delta_3+\Delta_4\wedge\omega_3-t\alpha\wedge\phi\;.
\]
The decompositions \eqref{eqn:W_i} of the components $W_3$ and $W_2^\pm$ correspond to projections
\begin{equation}
\label{eqn:Decompositions}\leftdoublebracket\Lambda^{2,1}_0\rightdoublebracket\to \omega_1^\perp\oplus
(\alpha\wedge\omega_1)^\perp\;, \quad [\Lambda^{1,1}_0]\to T\oplus \Lambda^2_-\;,
 \end{equation}
 respectively; the statement is now a straightforward consequence of \eqref{eqn:SU2Torsion}.
\end{proof}
\begin{remark}
Of the decompositions \eqref{eqn:Decompositions}, the first is not surjective; in other words, the components $\Xi_3$
and $\Delta_3$ are not independent.
\end{remark}

It is clear that one can write down a converse to Proposition~\ref{prop:QuotientTorsion}, because  the quotient
determines the intrinsic torsion of $M$; one can then characterize the $M$ with special intrinsic torsion in terms of
the quotient. For example, in the symplectic half-flat case one obtains this generalization of Lemma~\ref{lemma:GCY}:
\begin{proposition}\label{prop:GCYFromHypo}
$M$  is symplectic half-flat if and only if the quotient satisfies
\begin{equation}\label{eqn:GCYFromHypo}
\begin{gathered}
d\alpha=0\;,\quad d\omega_1=0\;,\quad  d\omega_2\wedge\alpha=t^2\omega_1\wedge\phi+2d\log t\wedge\omega_2\wedge\alpha\;,\\
d\omega_3=d\log t\wedge\omega_3-t\alpha\wedge\phi \;.
\end{gathered}
\end{equation}
\end{proposition}
\begin{proof}
Follows immediately from \eqref{eqn:HypoProduct2}.
\end{proof}
We now  consider the case where $M$ is integrable. In order to state our theorem, we need to introduce two differential
operators on  5-manifolds with an $\SU(2)$\nobreakdash-structure. The first one is $\partial_\alpha$, which maps a
function $f$ to $\langle \alpha,df\rangle$. Secondly, consider the endomorphism $J_3$ of $T^*N$ characterized by
\[J_3\alpha=0\;,\quad \omega_1\wedge\beta = \omega_2\wedge J_3\beta \text{ for } \beta\in\alpha^\perp\;;\]
we can then define an operator $d^c$ which maps a  function $f$ to  $d^c f=J_3df$.
\begin{theorem}
\label{thm:QuotientOfCY} If the $\SU(3)$-structure on $M$ is  integrable, the function $t$ is a solution of
\begin{equation}
\label{eqn:eqnfort}
 \partial_\alpha^2\log t -(\partial_\alpha \log t)^2 -2t^{-1}\norm{(d\log t)_{\Lambda^1}}^2=0\;.
\end{equation}
The intrinsic torsion is determined by $t$ as follows:  $\alpha$, $\omega_1$ and $\omega_2$ are closed, and $\omega_3$
satisfies
\begin{equation}
\label{eqn:eqnforomega3} d\omega_3=(d\log t)_{\Lambda^1}\wedge\omega_3 + \frac1{\partial_\alpha t}\,\alpha\wedge(2d\log
t\wedge d^c\log t-dd^c\log t)_{\Lambda^2_-}\;;
\end{equation}
moreover the ``curvature form'' is
 \begin{equation}
 \label{eqn:eqnforphi}
\phi = t^{-1}\partial_\alpha \log t\,\omega_3 - \frac1{t^2\partial_\alpha \log t}(2d\log t\wedge d^c\log t-dd^c\log
t)_{\Lambda^2_-} -2t^{-2}\alpha\wedge  d^c\log t\;.
\end{equation}
Conversely, let $N$ be a 5-manifold with an $\SU(2)$-structure $(\alpha,\omega_i)$ and a function $t$, where $\alpha$,
$\omega_1$ and $\omega_2$ are closed, and \eqref{eqn:eqnfort}, \eqref{eqn:eqnforomega3} are satisfied.  Then the
two-form $\phi$ defined by \eqref{eqn:eqnforphi} is closed; if the cohomology class $[\frac\phi{2\pi}]$ is an element
of $H^2(N,\Z)$, there is a circle bundle over $N$ on which  an integrable $\SU(3)$-structure is defined by
\eqref{eqn:HypoProduct2}, where $\eta$ is a connection form such that $d\eta=\phi$.
\end{theorem}
Locally, Theorem~\ref{thm:QuotientOfCY} is a characterization. Indeed, in the second part  one can restrict $N$ to a
contractible open subset $N'$, so that \[0=[\phi/2\pi]\in H^2(N',\Z)\;.\]
\begin{proof}
It is clear from  \eqref{eqn:quotient2} that $\alpha$, $\omega_1$ and $\omega_2$ are closed. Hence, the only
non-vanishing components of the intrinsic torsion are $\gamma_3$ and $\sigma_3^-$,  determined by
\[d\omega_3=\gamma_3\wedge\omega_3 + \alpha\wedge\sigma_3^-\;.\]
Using  \eqref{eqn:HypoProduct2}, we find
\[t^{-1}\omega_3\wedge (\gamma_3-d\log t) + \alpha\wedge(t^{-1}\sigma_3^-+\phi)=0\;.\]
Hence $\gamma_3 = (d\log t)_{\Lambda^1}$, and the component of $\phi$ in $\Lambda^2(\alpha^\perp)$ is determined by
\begin{align*}
\langle \phi,\omega_1\rangle &=0=\langle \phi,\omega_2\rangle \;,& \langle \phi,\omega_3\rangle&=2t^{-1}\partial_\alpha
\log t\;, & (\phi)_{\Lambda^2_-}&=-t^{-1}\sigma_3^-\;.
 \end{align*}
 It also follows from \eqref{eqn:HypoProduct2} that
\begin{align*}
\omega_1\wedge\phi+2t^{-3}dt\wedge\omega_2\wedge\alpha&=0\\
 \omega_2\wedge\phi-2t^{-3}dt\wedge\omega_1\wedge\alpha&=0
\end{align*}
 Therefore
\[
\alpha\hook\phi =-2t^{-2}d^c\log t\;.
\]
For brevity, we set $s=\partial_\alpha\log t$. By construction $\phi$ is closed, so
\begin{multline}
\label{eqn:Integrabledphi} 0=t^{-1}\bigl( -s\,d\log t\wedge\omega_3+ds\wedge \omega_3+ s\, (d\log t-s\,
\alpha)\wedge\omega_3 + s\,
\alpha\wedge\sigma_3^-+ \\
+d\log t\wedge\sigma_3^--d\sigma_3^-\bigr) +2t^{-2}\alpha\wedge(-2d\log t\wedge d^c\log t+dd^c\log t)\;.
\end{multline}
We can split \eqref{eqn:Integrabledphi}  into two equations by taking the wedge and the interior product with $\alpha$.
One of these is satisfied automatically: indeed, taking $d$ of $d\omega_3$ we find
\[0=\alpha\wedge(d\log t\wedge\sigma_3^- - d\sigma_3^- +ds\wedge\omega_3)\;,\]
so the right-hand side of \eqref{eqn:Integrabledphi} vanishes on wedging with $\alpha$. Taking the interior product
gives
\[(\partial_\alpha s -s^2)\omega_3 +2s\sigma_3^-+2t^{-1}(-2d\log t\wedge d^c\log t+dd^c\log t)_{\Lambda^2(\alpha^\perp)}\;.\]
It is now clear that $\sigma_3^-$ can be expressed in terms of $t$, giving \eqref{eqn:eqnforomega3}. Using the general
formula
\[\langle \beta\wedge J_3\beta,\omega_3\rangle = \norm{\beta}^2 - \langle \beta,\alpha\rangle^2\;,\]
we also deduce that $t$ satisfies \eqref{eqn:eqnfort}.

Conversely, suppose \eqref{eqn:eqnforomega3} and \eqref{eqn:eqnfort} are satisfied, and define  $\phi$ by
\eqref{eqn:eqnforphi}. The above calculations show that $\phi$ is closed and the construction of
Proposition~\ref{prop:Functor} defines an integrable $\SU(3)$-structure.
\end{proof}
\begin{remark}
The condition of Theorem~\ref{thm:QuotientOfCY} implies in particular  that 28 out of the 35 components of the
intrinsic torsion of $(N,\alpha,\omega_i)$ vanish. A similar construction was described in \cite{ApostolovSalamon},
starting with a 7-manifold with  holonomy $\Gtwo$, and defining an $\SU(3)$-structure  on the quotient. In that case,
the vanishing components of the intrinsic torsion of the quotient are also 28, though out of 42.
\end{remark}
In general \eqref{eqn:eqnfort} and \eqref{eqn:eqnforomega3} are not independent, because the norm on one-forms depends
on $\omega_3$. Motivated by this observation, we consider the  special case \[(d\log t)_{\Lambda^1}=0\;;\] in order to
apply Theorem~\ref{thm:QuotientOfCY}, we have to assume that the $\SU(2)$-structure is integrable. Let $x$ be a
coordinate in the direction of $\alpha$, so that $\alpha=dx$. Set $t=(1-x)^{-1}$; then \eqref{eqn:eqnfort} is
satisfied. Suppose that one has a circle bundle over $N$ with Chern class $[\frac1{2\pi}\omega_3]$; then the hypotheses
of Theorem~\ref{thm:QuotientOfCY} hold. Define a connection form $\eta$ such that $d\eta=\omega_3$; then
\[\omega=(1-x)\omega_3+\eta\wedge\alpha\;,\quad \Psi=(\omega_1+i\omega_2)\wedge\left(\eta+i(1-x)^2\alpha\right)\;,\]
defines an integrable $\SU(3)$-structure. One can actually prove that if the original 5\nobreakdash-manifold has
holonomy $\SU(2)$,
then the Calabi-Yau 6-manifold has  holonomy $\SU(3)$.\\

\bibliographystyle{plain}
\bibliography{contitomassini}

\begin{thebibliography}{10}

\bibitem{ApostolovSalamon}
V.~Apostolov and S.~Salamon.
\newblock {K\"a}hler reduction of metrics with holonomy {$G_2$}.
\newblock {\em Comm. Math. Phys.}, 246:43--61, 2004.

\bibitem{Bedulli}
L.~Bedulli.
\newblock {\em Tre-variet{\`a} di {C}alabi-{Y}au generalizzate}.
\newblock PhD thesis, Universit{\`a} di Firenze, 2004.

\bibitem{Blair}
D.E. Blair.
\newblock {\em {R}iemannian Geometry of Contact and Symplectic Manifolds}.
\newblock Birkh{\"a}user, 2002.

\bibitem{ChiossiFino}
S.~Chiossi and A.~Fino.
\newblock Conformally parallel {$G_2$} structures on a class of solvmanifolds.
\newblock To appear in {\em Math. Z.}

\bibitem{ChiossiSalamon}
S.~Chiossi and S.~Salamon.
\newblock The intrinsic torsion of {$SU(3)$} and {$G_2$ structures}.
\newblock In {\em Differential Geometry, Valencia 2001}, pages 115--133. World
  Scientific, 2002.

\bibitem{ChiossiSwann}
S.~Chiossi and A.~Swann.
\newblock {$G_2$}-structures with torsion from half-integrable nilmanifolds.
\newblock {\em J.Geom.Phys.}, 54:262--285, 2005.

\bibitem{ContiSalamon}
D.~Conti and S.~Salamon.
\newblock Generalized {K}illing spinors in dimension 5.
\newblock DG/0508375.

\bibitem{dBT3}
P.~de~Bartolomeis and A.~Tomassini.
\newblock On solvable generalized {C}alabi-{Y}au manifolds.
\newblock To appear in \emph{Ann. Inst. Fourier.}

\bibitem{dBT2}
P.~de~Bartolomeis and A.~Tomassini.
\newblock On the {M}aslov index of {L}agrangian submanifolds of generalized
  {C}alabi-{Y}au manifolds.
\newblock To appear in \emph{Int. J. of Math.}

\bibitem{FernandezGray}
M.~Fern{\'a}ndez and A.~Gray.
\newblock Riemannian manifolds with structure group {$\mathrm{G}_2$}.
\newblock {\em Annali di Mat. Pura Appl.}, 32:19--45, 1982.

\bibitem{Giovannini}
D.~Giovannini.
\newblock {\em Special structures and symplectic geometry}.
\newblock PhD thesis, Universit{\`{a}} degli studi di {T}orino, 2003.

\bibitem{HarveyLawson:CalibratedGeometries}
R.~Harvey and H.B. Lawson.
\newblock Calibrated geometries.
\newblock {\em Acta Math.}, 148:47--157, 1982.

\bibitem{Hitchin:TheGeometryOfThreeForms}
N.~Hitchin.
\newblock The geometry of three-forms in six dimensions.
\newblock {\em J. Differential Geom.}, 55:547--576, 2000.

\bibitem{Hitchin:StableForms}
N.~Hitchin.
\newblock Stable forms and special metrics.
\newblock In {\em Global Differential Geometry: The Mathematical Legacy of
  Alfred Gray}, volume 288 of {\em Contemp. Math.}, pages 70--89. American
  Math. Soc., 2001.

\end{thebibliography}
\end{document}